\documentclass[a4paper,11pt]{article}

\usepackage{amsfonts,amsthm,amsmath}

\usepackage{enumitem}

\usepackage[T1]{fontenc}

\newtheorem{teo}{Theorem}[section]
\newtheorem*{teo*}{Theorem}
\newtheorem{lem}[teo]{Lemma}
\newtheorem{cor}[teo]{Corollary}
\newtheorem{pro}[teo]{Proposition}
\newtheorem{conj}[teo]{Conjecture}

\theoremstyle{definition}

\theoremstyle{remark}
\newtheorem{exa}[teo]{Example}
\newtheorem{rem}[teo]{Remark}

\def\cB{\mathcal{B}}

\def\cH{\mathcal{H}}
\def\cI{\mathcal{I}}

\def\noi{\noindent}

\def\bdem{\begin{proof}}
\def\edem{\renewcommand{\qed}{\hfill $\blacksquare$}
\end{proof}}

\DeclareMathOperator{\ran}{ran}

\newcommand{\PI}[2]{\left\langle #1 , #2 \right\rangle}

\title{On a conjecture by Mbekhta about best approximation by  polar factors}
\author{Eduardo Chiumiento}
\date{}

\begin{document}

\maketitle

\begin{abstract}
The polar factor of a bounded operator acting on a Hilbert space is the unique partial isometry arising in the polar decomposition. 
It is well known that the polar factor might not be a best approximant to its associated operator in the set of all partial isometries, when the distance is measured in the operator norm. We show  that the polar factor of an arbitrary operator $T$ is a best approximant to $T$ in the set of all  partial isometries $X$ such that $\dim (\ker(X)\cap \ker(T)^\perp)\leq \dim (\ker(X)^\perp\cap \ker(T))$. We also provide a characterization of best approximations. This work is motivated by a recent conjecture by M. Mbekhta, which can be answered using our results.  
\end{abstract}

\bigskip




\section{Introduction}

 Let $\cH$ be a complex separable Hilbert space,  $\cB(\cH)$  the algebra of bounded linear operators  and  
 $\cI$  the set of all partial isometries on $\cH$.  The \textit{polar factor} of an operator $T$  is the unique $V \in \cI$ such $T=V|T|$ and $\ker(V)=\ker(T)$. Here, we write $|T|=(T^*T)^{1/2}$. The following conjecture was stated by M. Mbekhta  \cite{M}. 
\begin{conj}\label{op norm conj}
Let $T \in \cB(\cH)$ and $X_0 \in \cI$ such that $\ker(X_0)=\ker(T)$. The following conditions are equivalent:
\begin{itemize}
\item[(i)] $X_0$ is the polar factor of $T$;
\item[(ii)] $\|T-X_0\|=\min\{    \|T- X \| \, : \, X \in \cI,  \, \ker(X)=\ker(T) \}.$
\end{itemize}
\end{conj}
The norm considered is the usual operator norm (or spectral norm). M. Mbekhta proved as a partial result that the first item implies the second when $T$ is injective.  In this work we show that this implication holds true for an arbitrary operator $T$. We also show that the 
converse implication is false.  

We actually prove in Section \ref{best a p} that the polar factor is a best approximant to its associated operator in a larger set of partial isometries. Instead of fixing the kernel of the partial isometries, our set is given in terms of the dimension of certain subspaces. More precisely, given $P,Q$ two orthogonal projections on  $\cH$, we set
$$
j(P,Q):=\dim(\ran(P) \cap \ker(Q)) - \dim(\ker(P)  \cap \ran(Q)) \in [-\infty,\infty], 
$$
if one of these dimensions is finite, and $j(P,Q)=0$ if both dimensions are infinite. 
In Theorem \ref{principal thm} we prove the following best approximation property of the polar factor $V$ of an operator $T$: 
\begin{equation}\label{ref min}
\| T- V \|  = \min \{   \|T-X \| \,: \, X \in \cI, \, j(V^*V, X^*X) \leq 0 \}.
\end{equation}
From this result, we   obtain that the first item implies the second in Conjecture \ref{op norm conj}. Indeed, note that $j(V^*V,X^*X)=0$ if $\ker(X)=\ker(T)$ $(=\ker(V))$. 
Then we find a necessary and sufficient condition for the polar factor to be a best approximant to its associated operator in the set of all partial isometries. The reduced minimum modulus $\gamma(T)$ of the operator $T$ naturally shows up as in this condition. More precisely, we prove  in Proposition \ref{spectral cond} that   $\|T\|-1<1- \gamma(T)$, or $\gamma(T)<1/2$, if and only if the polar factor $V$ is a best approximant to $T$ in the set $\cI$. If our condition fails, then other partial isometries, which are best approximants in $\cI$ and different from the polar factor, can be explicitly    constructed (see Remark \ref{other minimizers}).

We have already mentioned that the converse of Conjecture \ref{op norm conj} does not hold true; a counterexample  is given at the beginning of Section \ref{charact}. 
 This  motivates the characterization  we give in Theorem \ref{caract thm} of those  partial isometries attaining the minimum in (\ref{ref min}). It relies on characterizations of operators satisfying the equality case in the triangle inequality \cite{BB02, L99} combined with our previous results. 

We end this section with  remarks and connections to previous works. It is worth pointing  out that the definition of  $j(P,Q)$ already appeared in the literature. If the operator $QP|_{\ran(P)}:\ran(P) \to \ran(Q)$ is Fredholm, then $j(P,Q)$ coincides with its Fredholm index. In this case, the pair $(P,Q)$ is called a Fredholm pair, a notion studied in  \cite{AS94, ASS94}. A fundamental result in this regard, which applies also when $(P,Q)$ is not necessary Fredholm, is that the condition $j(P,Q)=0$ is  equivalent  to unitary equivalence of the projections with a unitary permuting them (see \cite{S17, WDD09}). 


On the false implication of the conjecture, we remark that the operator norm is neither strictly convex nor differentiable. These facts
are essential to obtain information about the minimizers. For instance, this implication of the conjecture holds true when the operator norm is replaced by the Hilbert-Schmidt norm for arbitrary Hilbert spaces \cite{C19, FPT02}, or by any strictly convex unitarily invariant norm for finite dimensional Hilbert spaces \cite{AC18}.  In these works, best approximation by partial isometries with the Hilbert-Schmidt norm has deserved special attention due its importance for frame theory. In particular we observe that in \cite{C19} previous results using the notion of Fredholm pairs of projections were given; though for dealing with the Hilbert-Schmidt norm the proofs depend on different techniques like majorization of singular values, which do not apply to the present work.  Best approximation by partial isometries with the operator norm was studied by   P.W. Wu \cite{W86}, where among other results, it was shown that the distance of an operator to the set of all partial isometries is always attained.  However, the role of the polar factor in best approximation by partial isometries was not considered in that work.  Other related articles concern with best approximation by partial isometries in the matricial case and its applications \cite{LZ06}, best approximation using the Schatten $p$-norms of operators \cite{M89}, and the polar decomposition of products of projections \cite{CM11}.




\section{Best approximation by the polar factor}\label{best a p}
An operator $X \in \cB(\cH)$ is called a partial isometry if $XX^*X=X$, or equivalently, if $X^*X$ is an orthogonal projection. This is 
also equivalent to have $\|X\xi\|=\|\xi\|$, for every vector $\xi \in \ker(X)^\perp$. 
As we mentioned above,  $\cI$ denotes the set consisting of all partial isometries. We observe that $X \in \cI$ if and only if $X^* \in \cI$. 
Lastly, we note for later use that when $X \in \cI$, the projection $X^*X$, called the initial projection, projects onto $\ker(X)^\perp$; meanwhile $XX^*$, the final projection, projects onto $\ran(X)$.

Let $\sigma(T)$ denote the spectrum of an operator $T$. Recall the definition of the reduced minimum  modulus $\gamma(T)$ of an operator $T\neq 0$  (see, e.g. \cite{A85}):
\begin{align}
\gamma(T) & =\inf\{ \|T\xi \| \, : \, \xi \in \ker(T)^\perp, \, \|\xi\|=1  \} \nonumber \\
& = \inf \sigma(|T|)\setminus \{ 0 \}  \label{red mod min}.
\end{align}

\begin{lem}\label{pf distance}
Let $T \in \cB(\cH)$, $T\neq 0$, with polar decomposition $T=V|T|$. Then
$$
\|T - V \|= \max \{ 1- \gamma(T), \, \|T\| - 1\}. 
$$
\end{lem}
\begin{proof}
We write $P=V^*V$ for the initial projection of the partial isometry $V$. Since $\overline{\ran}(|T|)=\ker(|T|)^\perp=\ker(V)^\perp=\ran(P)$, the partial isometry $V$ acts isometrically on $\ran(|T|)$, and it follows that
$\|T- V\|=\| V(|T| - P)\|=\| |T| - P \|$. Noting that $\ker(|T|-P)^\perp \subseteq \ran(P)$ and $\ran(|T| - P)\subseteq \ran(P)$, we have $\| |T| - P \|=\| (|T| - P)|_{\ran(P)} \|$, where the second norm is taken in the invariant subspace $\ran(P)$. If we put  $T_0:= |T||_{\ran(P)}: \ran(P) \to \ran(P)$, then $\| |T|-P\|=\|T_0-I\|$, where $I$ is the identity on $\ran(P)\neq \{ 0 \}$. Therefore an application of the spectral theorem on $\ran(P)$ gives
\begin{align*}
\||T| - P \| & = \| T_0 - I \|\\
& =\sup \{  |\lambda - 1| \, : \lambda \in \sigma(T_0) \} \\
&  = \max \{  1- \inf \sigma(T_0), \, \sup \sigma(T_0) -1  \} \\
& =\max \{  1- \gamma(T), \, \|T\| - 1  \}.
\end{align*}
In the last equality we have used that $\sup \sigma(T_0)=\|T\|$ and $\inf \sigma(T_0)=\gamma(T)$. The proofs of these facts are straightforward;  one only needs to note that 
$\sigma(T_0)\setminus \{ 0\}=\sigma(|T|)\setminus \{ 0\}$ because $\ran(P)$ is an invariant subspace of $|T|$ and $\ran(P)^\perp=\ker(|T|)$.
\end{proof}


\begin{lem}{\cite[Thm. 3.4]{W86}}\label{Wu thm}
Given  $T, S \in \cB(\cH)$, then
$$
\|T-S\|\geq \sup_{\lambda \in \sigma(|T|)} \inf_{\mu \in \sigma(|S|)} \{ \, \lambda, \, |\lambda - \mu| \, \}.
$$
Moreover, if $\dim \ker(S) \geq \dim \ran(S)^\perp$, then
$$
\|T-S\|\geq \sup_{\lambda \in \sigma(|T|)} \inf_{\mu \in \sigma(|S|)}  |\lambda - \mu| .
$$
\end{lem}
\begin{rem}
The dimensions of the subspaces $\ker(S)$ and $\ran(S)^\perp$ in the second part of the statement could be infinite.
\end{rem}



Our main result on best approximation by the polar factor is the following.

\begin{teo}\label{principal thm}
Let $T \in \cB(\cH)$ with polar decomposition $T=V|T|$. Then
\begin{align*}
\| T- V \| & = \min \{   \|T-X \| \,: \, X \in \cI, \, j(V^*V, X^*X) \leq 0 \} \\
& = \min \{   \|T-X \| \,: \, X \in \cI, \, j(VV^*, XX^*) \leq 0 \}. 
\end{align*}
\end{teo}
\begin{proof}
 It is not difficult to check that $T^*=V^*|T^*|$ is the polar decomposition of $T^*$.  Thus, the second characterization using the final projections can be derived from the first characterization using the initial projections. Also we observe that the result clearly follows when $T=0$. In order to prove the first characterization for $T\neq 0$, we divide the proof into two cases.  In the first case we suppose that $\|T-V\|=\|T\| - 1$ in Lemma \ref{pf distance}. In particular, this implies $\|T\|\geq 1$.  For any partial isometry $X$ we have that $X^*X$ is a projection, so that $\sigma(|X|)=\{0, \, 1 \}$.  According to Lemma   \ref{Wu thm}, and noting that $\|T\| \in \sigma(|T|)$, we have
\begin{align*}
\|T- X\| & \geq \sup_{\lambda \in \sigma(|T|)} \inf_{\mu \in \{ 0, \, 1 \}} \{  \lambda, \, |\lambda - \mu|  \} \\
& \geq \inf_{\mu \in \{ 0, \, 1 \}} \{  \|T\|, \, \|T\| - \mu  \} \\
& = \|T\| - 1 = \|T -V\|.
\end{align*}
This proves the first case. 

In the second case we assume that $\|T-V\|=1 - \gamma(T)$ in Lemma \ref{pf distance}. Thus, $\gamma(T) \in [0,1]$.  Pick a partial isometry $X$
such that $j(V^*V, X^*X) \leq 0$. For notational simplicity, we write $P=V^*V$ and $Q=X^*X$. We divide this part of the proof in three steps.


\smallskip

Step 1. If we suppose that there is a vector $\xi \in \ker(P)\cap \ran(Q)$,  $\|\xi\|=1$, then
\begin{align}
\|T- X\| & \geq \| (T-X)\xi\|=\|X \xi\| =1 \nonumber \\
&  \geq 1 - \gamma(T)=\|T-V\|. \label{1case}
\end{align}
Thus, we can assume that $\ker(P)\cap \ran(Q)=\{0\}$ for the rest of the proof. By the condition $j(P,Q)\leq 0$, we also have $\ker(Q) \cap \ran(P)=\{ 0\}$. Hence,  $G:=PQ|_{\ran(Q)}:\ran(Q) \to \ran(P)$ is an  injective operator with dense range. 

\smallskip

Step 2. Next we suppose that $\overline{\ran}(T)\cap \ran(X)^\perp \neq \{ 0\}$. We are going to prove the desired conclusion under this assumption. 
Take a sequence $(\xi_n)$ in $\ker(T)^\perp$, $T\xi_n \to \eta \in \ran(X)^\perp$, $\eta \neq 0$.  Since the range of the operator $G$ is dense in $\ran(P)=\ker(T)^\perp$, then there is a sequence $(\zeta_n)$ in $\ran(Q)$ such that $\| G \zeta_n - \xi_n\| \leq 1/n$.  
Observe that
\begin{align*}
\| \zeta_n\|^2 \|T-X\|^2 & \geq \| \zeta_n \|^2 \|TPQ - XQ \|^2 \geq \|T G \zeta_n  - X \zeta_n \|^2 \\
& = \|T G \zeta_n \|^2 + \|X \zeta_n \|^2 - 2\Re\PI{TG\zeta_n}{X\zeta_n}, 
\end{align*}
which gives
\begin{align}
\| T - X \|^2 & \geq \|TPQ - X\|^2 \nonumber \\
& \geq \left\| TG \frac{\zeta_n}{\|\zeta_n\|} \right\|^2 + \left\|X \frac{\zeta_n}{\|\zeta_n\|} \right\|^2 - \frac{2}{\|\zeta_n\|^2} \Re\PI{TG\zeta_n}{X\zeta_n} \nonumber\\
& \geq 1 - \frac{2}{\|\zeta_n\|} \, \Re \PI{TG\zeta_n}{X \frac{\zeta_n}{\|\zeta_n\|}}. \label{number1}
\end{align}
In the last inequality we have used that $\|\zeta_n \|^{-1}\zeta_n$  are unit vectors in $\ran(Q)=\ker(X)^\perp$. We have to show that 
the second term in (\ref{number1}) goes to zero. To this end, we first note that $\|T\xi_n - \eta\| \geq \| \eta \| - \|T\|\|\xi_n\|$, which implies
$\|\xi_n\| \geq \|T\|^{-1}(\|\eta\| - \|T\xi_n - \eta\|)$.   Using that $1/n \geq \|G\zeta_n - \xi_n\|\geq \|\xi_n \| - \|G\|\|\zeta_n\|$, 
and noting that $\|G\|\leq 1$, we get
$$
\| \zeta_n\| \geq \|\xi_n\| - \frac{1}{n} \geq  \frac{\|\eta\|}{\|T\|} - \frac{\|T\xi_n - \eta\|}{\|T\|} - \frac{1}{n}.
$$
Thus, for $n\geq 1$ large enough, it follows that $\|\zeta_n\| \geq C$, for some constant $C>0$. Now we use that $\eta \in \ran(X)^\perp$ to compute the limit:
\begin{align*}
\left|\|\zeta_n\|^{-1} \PI{TG\zeta_n}{X \frac{\zeta_n}{\|\zeta_n\|}}\right| & \leq  C^{-1}\left(\left| \PI{TG\zeta_n - \eta}{X\frac{\zeta_n}{\|\zeta_n\|}}\right| + \left|\PI{\eta}{X\frac{\zeta_n}{\|\zeta_n\|}}\right| \right)\\
& \leq C^{-1}\|TG \zeta_n - \eta \| \, \left\|X\frac{\zeta_n}{\|\zeta_n\|}\right\| \\
& \leq C^{-1} (\| T \| \,\|G\zeta_n - \xi_n \| + \|T\xi_n - \eta\|) \to 0.
\end{align*}
Taking  limit in (\ref{number1}) we find that
$\|T- X\| \geq 1 \geq 1 -\gamma(T)=\|T-V\|$. 

\smallskip

Step 3. We denote by  $F$ and $E$ the orthogonal projections onto $\overline{\ran}(T)$ and $\ran(X)$, respectively.  From the previous steps, we can assume that  $G=PQ|_{\ran(Q)}:\ran(Q) \to \ran(P)$ is injective with dense range, and the operator $R:=EF|_{\ran(F)}:\ran(F)\to \ran(E)$ is injective. We consider the restrictions $T_0:=ET|_{\ran(Q)}:\ran(Q) \to \ran(E)$ and $X_0:=EX|_{\ran(Q)}:\ran(Q)\to \ran(E)$. These  clearly satisfy $\| E(T-X)Q\|=\| T_0 - X_0\|$, where the first operator norm is taken as an operator on $\cH$, and the second as an operator from $\ran(Q)$ to $\ran(E)$. Before the forthcoming inequalities, it is also convenient to  observe two facts. First, 
the operator $X_0$ is an isometric isomorphism, which implies that  $\sigma(|X_0|)=\{ 1\}$. Second, we note that $\gamma(T_0) \in \sigma(|T_0|)$ by the characterization  (\ref{red mod min}) of the reduced minimum  modulus.  Then an application of the second statement in Lemma \ref{Wu thm} yields
\begin{align}
\|T-X\| & \geq \| T_0 - X_0\| \nonumber \\
& \geq \sup_{\lambda \in \sigma(|T_0|)} \inf_{\mu \in \sigma(|X_0|)}  |\mu - \lambda| \nonumber  \\
& = \sup_{\lambda \in \sigma(|T_0|)} |1- \lambda| \nonumber \\
& \geq  1 - \gamma(T_0)  \nonumber \\
& \geq 1 - \gamma(T) = \|T- V\|. \label{number3}
\end{align}
It remains to be shown the last inequality used above: $\gamma(T_0)\leq \gamma(T)$. To see this fact, note that $T_0$ is injective. This follows immediately using that $T_0=ET|_{\ran(Q)}=EFTPQ|_{\ran(Q)}=RTG|_{\ran(Q)}$, and noting that $R$ and $G$ are injective with $\overline{\ran}(G)=\ran(P)=\ker(T)^\perp$. Since $\|R\|\leq 1$ and $\|G\|\leq 1$, then
\begin{align}
\gamma(T_0) & = \inf \{ \|RTG\zeta \| :  \zeta \in \ran(Q), \, \|\zeta \|=1 \} \nonumber \\ 
& \leq \inf \{ \|T G \zeta \| :  \zeta \in \ran(Q), \, \| \zeta\|=1 \} \nonumber \\
& = \inf \left\{ \left\|T \frac{G \zeta}{\|G\zeta\|} \right\| \|G\zeta \| :  \zeta \in \ran(Q), \, \| \zeta\|=1 \right\} \nonumber  \\
& \leq \inf \left\{ \left\|T \frac{G \zeta}{\|G\zeta\|} \right\| :  \zeta \in \ran(Q), \, \| \zeta\|=1 \right\} = \gamma(T), \label{number4}
\end{align}
where the last equality follows by using again the fact  $\overline{\ran}(G)=\ker(T)^\perp$. This finishes the proof.
\end{proof}


Now we can prove that the first item implies the second in Conjecture \ref{op norm conj}. 

\begin{cor}
Let $T \in \cB(\cH)$ with polar decomposition $T=V|T|$. Then
$$
\|T-V\|=\min \{  \|T- X\| \, : \, X \in \cI, \, \ker(X)=\ker(T) \}.
$$
\end{cor}
\begin{proof}
Since $\ker(X)=\ker(T)$ $(=\ker(V))$,  then $j(V^*V,X^*X)=0$. Then the result follows immediately from Theorem \ref{principal thm}.
\end{proof}


The distance of an arbitrary operator to the set of all partial isometries is always attained (see \cite[Thm. 3.6]{W86}).  
Based on this result, we can give necessary and sufficient conditions on the spectrum of an operator to guarantee that its polar factor becomes a  best approximant in the set of all partial isometries.

\begin{pro}\label{spectral cond}
Let $T \in \cB(\cH)$, $T\neq 0$,  with polar decomposition $T=V|T|$. The following conditions are equivalent:
\begin{itemize}
\item[(i)] $\|T\|-1 \geq 1- \gamma(T)$, or  $\gamma(T)\geq 1/2$;
\item[(ii)] $\|T-V\|=\min \{   \| T- X\| \, : \,  X \in \cI  \}$.
\end{itemize}
\end{pro}
\begin{proof}
Consider the function $f(\lambda)=\min \{ \lambda, \, |1-\lambda| \}$, and put $i_0:=\inf\{ \|T-X\| \, : \, X \in \cI\}$. In \cite[Thm. 3.6]{W86} P.Y. Wu proved that this infimum is attained. Moreover, it can be computed in terms of the spectrum of $|T|$ as $i_0=\sup_{\lambda \in \sigma(|T|)} f(\lambda)$. 
As in the proof of Lemma \ref{pf distance}, we set $P=V^*V$ and $T_0=|T||_{\ran(P)}:\ran(P)\to \ran(P)$. Recall that in the proof of the aforementioned lemma, we have shown that $\|T-V\|=\sup_{\lambda \in \sigma(T_0)}|1- \lambda|$, $\sigma(T_0)\setminus\{0\}=\sigma(|T|)\setminus\{0\}$ and $\inf \sigma(T_0)=\gamma(T)$. Under the assumptions $\|T\|-1 \geq 1- \gamma(T)$, or  $\gamma(T)\geq 1/2$, it follows that 
$$
i_0=\sup_{\lambda \in \sigma(|T|)} f(\lambda)= \sup_{\lambda \in \sigma(T_0)} |1- \lambda|=\|T-V\|.
$$
This proves  one implication. 

In order to prove the converse, we assume that $\gamma(T)<1/2$ and $\|T\|-1< 1 - \gamma(T)$. By Lemma \ref{pf distance}, $\|T-V\|=1-\gamma(T)$. Note that the function $f$ can attain its maximum restricted to $\sigma(|T|)$ in any of the following intervals $[0,1/2]$, $[1/2,1]$ and $[1,\infty)$. Then we have that $i_0=\sup_{\lambda \in \sigma(|T|)} f(\lambda)=\max\{ a, \, 1-b, \, \|T\| -1\}$ for some positive numbers $a \leq 1/2$ and $b\geq 1/2$. We consider the three cases. In the first case,  $a<1-\gamma(T)$ by the assumption $\gamma(T)<1/2$. Similarly, in the second case $1-b<1-\gamma(T)$.
The last case uses the assumption $\|T\|-1< 1 - \gamma(T)$. Hence $i_0<\|T-V\|$.
\end{proof}

\begin{rem}\label{other minimizers}
In the case where  $\gamma(T)<1/2$ and $\|T\|-1< 1 - \gamma(T)$, we have shown in the above proof that $i_0=\inf\{ \|T-X\| \, : \, X \in \cI\}<\|T-V\|$. The infimum is attained, so this means that there exists a partial isometry $X_0$ such that $i_0=\|T-X_0\|<\|T-V\|$. According to Theorem \ref{principal thm}, it must be $j(V^*V, X_0^*X_0)>0$. It is interesting to recall how $X_0$ is constructed in \cite[Thm. 3.6]{W86}. 
For this, consider the function  $\phi(t)=1 - \chi_{(0,1/2)}(t)$. Then, $X_0$ is defined using Borel functional calculus by $X_0=V\phi(|T|)$.
Note that by the condition $\gamma(T)<1/2$, $\phi(|T|)$ turns out to be a proper subprojection of $V^*V$.  
\end{rem}

\section{Characterization of best approximations}\label{charact}

The following example shows that the second item does not imply the first item in  Conjecture \ref{op norm conj}.

\begin{exa}\label{ex1}
Consider the following matrices 
$$
T=\begin{pmatrix}
a &  0   &  0  \\  0 &  1  & 0 \\ 0   &    0   &  1
\end{pmatrix},
\, \, \, \,  a > 3;
\, \, \, \, \, \, \, \, \, \,  
X_0=\begin{pmatrix} 
1 &  0   &  0  \\  0 &  0  & -1 \\ 0   &    -1   &  0
\end{pmatrix}. 
$$
The polar factor of $T$ is the identity matrix $I$, and $\|T - I\|=a-1$. Clearly, $\ker(T)=\ker(X_0)=\{ 0\}$ and
\begin{align*}
\|T- X_0\| & =\max \left\{ \, a-1 \, , \, \left\|\begin{pmatrix} 1  & 1  \\  1  &   1 \end{pmatrix}\right\| \, \right\} \\
& =\max\{ \, a - 1 \, , \, 2 \, \}=a-1.
\end{align*}
Thus, we have 
$$\|T- X_0\|=\|T - I \|=\min \{   \| T- X \| \, : \, X \in \cI, \, \ker(X)=\{ 0\}  \},$$ 
and $X_0 \neq I$.
\end{exa}

It is then natural to study those partial isometries whose distance to a fixed operator coincides with the distance of the operator to its polar factor. That is, those partial isometries that attain the minimum in Conjecture \ref{op norm conj}, or more generally, the minimum in Theorem \ref{principal thm}. In this direction we have the following result.

\begin{teo}\label{caract thm}
Let $T \in \cB(\cH)$ with polar decomposition $T=V|T|$. Let $X_0$ be a partial isometry satisfying $j(V^*V,X_0^*X_0)\leq 0$. Then, 
\begin{equation}\label{attains min}
\|T-X_0\|=\min \{   \|T-X \| \,: \, X \in \cI, \, j(V^*V, X^*X) \leq 0 \}
\end{equation}
if and only if any of the following conditions hold:
\begin{itemize}
\item[(i)] There exist unit vectors $(\xi_n)$ such that $\|X_0\xi_n\|\to 1$ and  $\|T-X_0\|X_0\xi_n - (T-X_0)\xi_n \to 0$.
\item[(ii)] There exist unit vectors $(\xi_n)$ such that $\|X_0\xi_n\|\to 1$,  $\|T-X_0\|X_0\xi_n + (T-X_0)\xi_n \to 0$ and $(\gamma(T)X_0-T)\xi_n \to 0$.
\end{itemize}
\end{teo}
\begin{proof}
Suppose that $X_0$ is a partial isometry, $j(V^*V,X_0^*X_0)\leq0$, and $X_0$ attains the minimum in (\ref{attains min}). Then, $\|T-X_0\|=\|T-V\|$ by Theorem \ref{principal thm}. According to Lemma \ref{pf distance}, we have to consider the cases $\|T-V\|=\|T\| -1$ and $\|T-V\|=1- \gamma(T)$. In the first case,  we have $\|T-X_0\|+\|X_0\|=\|T\|$. That is, equality in the triangle inequality for the operator norm. From \cite[Thm. 1]{L99} we have that there exists a sequence of unit vectors $(\xi_n)$ such that $\|X_0-T\|X_0\xi_n -(T-X_0)\xi_n \to 0$. An examination of the proof of the quoted result yields that the sequence  $(\xi_n)$ satisfies $\|(T-X_0) \xi_n\|\to \|T-X_0\|$. From this fact, it follows that $\|X_0\xi_n\| \to 1$. 

In the second case, we assume that  $\|T-X_0\|=\|T-V\|=1 -\gamma(T)>\|T\|-1$. As in the proof of Theorem \ref{principal thm}, we have to consider three steps. Put $P=V^*V$ and $Q=X_0^*X_0$. 
\smallskip

Step 1. Suppose that $\ker(P)\cap \ran(Q)\neq \{ 0\}$. From the inequalities 
(\ref{1case}) and the assumption $\|T-X_0\|=\|T-V\|$, it must be $\gamma(T)=0$ and $\|T-X_0\|=1$. Then, pick a unit vector $\xi \in \ker(P)\cap \ran(Q)$, and note that $\|T-X_0\|X_0\xi + (T-X_0)\xi=0$.  This is the required conclusion. For the remainder of the proof of this item we will suppose that
$\ran(Q) \cap \ker(P)=\{ 0\}$. Since $j(P,Q)\leq 0$, then also $\ran(P)\cap \ker(Q)=\{ 0\}$.

\smallskip

Step 2. We further assume that $\overline{\ran}(T)\cap \ran(X_0)^\perp \neq \{ 0\}$.  At the end of the proof of the second step in Theorem \ref{principal thm} we have shown that $\|T-X_0\| \geq 1 \geq 1-\gamma(T)=\|T-V\|$. Again, since $\|T-X_0\|=\|T-V\|$, then $\gamma(T)=0$ and $\|T-X_0\|=1$.
From the inequalities (\ref{number1}), there exist unit vectors $\xi_n:=\|\zeta_n\|^{-1}\zeta_n \in \ran(Q)$  such that 
$$
1 = \| T-X_0\|^2 \geq \|(T-X_0)\xi_n\|^2 =  \|T \xi_n\|^2 + \|X_0\xi_n\|^2 + A_n \geq  1 + A_n,   
$$
where $(A_n)$  is a sequence converging to zero. Hence $(\xi_n)$ are unit vectors in $\ran(Q)$ satisfying $\|(T-X_0)\xi_n\| \to \|T-X_0\|$ and $\|T\xi_n\|\to \gamma(T)$. We will see that this condition is sufficient to finish the proof.

\smallskip

Step 3. Now we are under the same assumptions and notation of the third step in the proof of Theorem \ref{principal thm}.
Using  that $\|T-X_0\|=\|T-V\|$ in the inequalities (\ref{number3}), it follows  that $\gamma(T_0)=\gamma(T)$. Therefore all the inequalities in  (\ref{number4}) become equalities, which in particular implies that there are unit vectors $(\xi_n)$ in $\ran(Q)$ such that 
$\|T \xi_n\|=\|TG\xi_n\| \to \gamma(T)$. Furthermore, $\|T-X_0\|\geq \|(T-X_0)\xi_n\| \geq 1 - \|T\xi_n\|$ and $\|T-X_0\|=1- \gamma(T)$, implies
$\|(T-X_0)\xi_n\| \to \|T-X_0\|$.

\smallskip

We have seen in both step 2 and 3  that there are unit vectors $(\xi_n) \subseteq \ran(Q)$ such that $\|T\xi_n\| \to \gamma(T)$ and $\|(T-X_0)\xi_n \|\to \|T-X_0\|$. Now note that 
$$
\|(\gamma(T)X_0 + T)\xi_n\| \geq | \, \| (\gamma(T)+1)X_0\xi_n \| -  \| (X_0-T)\xi_n \| \, |\to |1+ \gamma(T)- (1-\gamma(T)) |=2\gamma(T).
$$ 
Therefore, by the parallelogram law, 
\begin{align*}
\|(\gamma(T)X_0-T)\xi_n\|^2=2(\gamma(T)^2 + \|T\xi_n\|^2)- \|(\gamma(T)X_0 + T)\xi_n\|^2 \to 0. 
\end{align*}
Now the  remaining condition follows:
\begin{align*}
\| \, (\|T-X_0\|X_0 + (T-X_0))\xi_n \| & \leq   \|(\|T-X_0\| + (\gamma(T) -1 ) )X_0\xi_n\|     +  \|(T- \gamma(T)X_0)\xi_n\| \\
& = \|(T- \gamma(T)X_0)\xi_n\| \to 0.
\end{align*}

To prove the converse, take $(\xi_n)$ unit vectors satisfying the conditions of item $(i)$.  Observe that
$$
\PI{(T-X_0)\xi_n}{X_0\xi_n}=\PI{((T-X_0) - \|T-X_0\|X_0)\xi_n}{X_0\xi_n} + \|T-X_0\|\|X_0\xi_n\|^2 \to \|T-X_0\|.
$$
In \cite[Thm. 2.1]{BB02} the authors  proved that given two operators $A,B$, then $\|A+B\|=\|A\|+\|B\|$ if and only if $\|A\|\|B\| \in \overline{W(A^*B)}$.
Here $W(C)=\{ \PI{C\xi}{\xi} \, : \, \|\xi\|=1\}$ is the numerical range of an operator $C$. The above computation means that 
$\|X_0\|\|T-X_0\|=\|T-X_0\| \in \overline{W(X_0^*(T-X_0))}$. Thus, we have $\|T-X_0\| + \|X_0\|=\|T\|$. Then, $\|T-X_0\| =\|T\|-1 \leq \|T-V\|$ by Lemma \ref{pf distance}. Since we are assuming that $j(V^*V,X_0^*X_0)\leq 0$, Theorem \ref{principal thm} implies that $\|T-X_0\|=\|T-V\|$, and thus $X_0$ attains the minimum in (\ref{attains min}).

Now suppose that $(\xi_n)$ are unit vectors satisfying the conditions in item $(ii)$. Note that
\begin{align*}
\|(T-X_0)\xi_n\| \leq \|(T- \gamma(T)X_0)\xi_n\| + \|(\gamma(T)-1)X_0\xi_n\| \to |1- \gamma(T)|.
\end{align*}
Since $\|(T-X_0)\xi_n\|\to \|T-X_0\|$, then $\|T-X_0\|\leq |1-\gamma(T)|\leq \|T-V\|$. Again this implies that $\|T-X_0\|=\|T-V\|$ and $X_0$ attains the minimum.
\end{proof}

\begin{rem}
Both of the conditions $\|T-X_0\|X_0\xi_n + (T-X_0)\xi_n \to 0$ and $(\gamma(T)X_0-T)\xi_n \to 0$  in item $(ii)$ of the above theorem  are used to prove that $X_0$ is a minimizer. We remark that only one of these conditions is not sufficient to obtain the same conclusion.
 For instance, take 
$$
T=\begin{pmatrix}
1 &  0     \\  0 &  1/2 
\end{pmatrix},
\, \, \, \,  
\, \, \, \, \, \, \, \, \, \,  
X_0=\begin{pmatrix} 
-1 &  0    \\  0  & 1 
\end{pmatrix}. 
$$
Then, $\|T-V\|=1/2$, $\|T-X_0\|=2$ and $\|T-X_0\|X_0 \xi_1 + (T-X_0) \xi_1=0$, $\xi_1=(1,0)$. Also note that $(\gamma(T)X_0 - T)\xi_2=0$,
$\xi_2=(0,1)$.
\end{rem}



\subsection*{Acknowledgment}
I thank Mostafa Mbekhta for sharing his work \cite{M} with me before publication and valuable comments. Also I would like to thank 
Gustavo Corach for several conversations on this subject. Finally I thank the referees for their constructive comments. This research was supported by Grants CONICET (PIP 2016 0525), ANPCyT (2015 1505/ 2017 0883) and FCE-UNLP (11X829).

\bigskip

{\sc (Eduardo Chiumiento)} {Departamento de  Matem\'atica \& Centro de Matem\'atica La Plata, FCE-UNLP, Calles 50 y 115, 
(1900) La Plata, Argentina  and Instituto Argentino de Matem\'atica, `Alberto P. Calder\'on', CONICET, Saavedra 15 3er. piso,
(1083) Buenos Aires, Argentina.}     
                                               
\noi e-mail: {\sf eduardo@mate.unlp.edu.ar}


\begin{thebibliography}{XX}

\bibitem{AS94} W.O. Amrein, K.B. Sinha, {\sc On pairs of projections in a Hilbert space}, Linear Algebra Appl. 208/209 (1994). 425-435.


\bibitem{AC18} J. Antezana, E. Chiumiento, {\sc Approximation by partial isometries and symmetric approximation of finite frames}, J. Fourier Anal. Appl. 24 (2018), no. 4, 1098-1118. 


\bibitem{A85} C. Apostol, {\sc The reduced minimum modulus}, Michigan Math. J. 32 (1985), 279-294.


\bibitem{ASS94} J. Avron, R. Seiler, B. Simon, {\sc The index of a pair of projections}, J. Funct. Anal. 120 (1994), no. 1, 220-237. 


\bibitem{BB02} M. Barraa, M. Boumazgour, {\sc Inner derivations and norm equality}, Proc. Amer. Math. Soc. 130 (2002), no. 2,  471-476.


\bibitem{C19} E. Chiumiento, {\sc Global symmetric approximation of frames},  J. Fourier Anal. Appl. 25  (2019), 1395-1423. 


\bibitem{CM11} G. Corach, A. Maestripieri, {\sc Products of orthogonal projections and polar decompositions}, Linear Algebra Appl. 434 (2011), 1594-1609. 

\bibitem{FPT02} M. Frank, V. Paulsen, T. Tiballi, {\sc Symmetric Approximation of frames and bases in Hilbert Spaces}, Trans. Amer. Math. Soc. 354 (2002), 777-793.

\bibitem{LZ06} B. Laszkiewicz, K. Zi\c{e}tak, {\sc Approximation of matrices and a family of Gander methods for polar decomposition}, BIT 46 (2006), no. 2, 345-366.



\bibitem{L99} C.S. Lin, {\sc The unilateral shift and a norm equality for bounded linear operators}, Proc. Amer. Math.
Soc. 127  (1999), no. 6, 1693-1696.

\bibitem{M89} P.J. Maher, {\sc Partially isometric approximation of positive operators}, Illinois J. Math. 33
(1989),  227-243.



\bibitem{M}  M. Mbekhta, {\sc Approximation of the polar factor of an operator acting on a Hilbert space}, J. Math. Anal. Appl. 487 (2020), 123954. 

\bibitem{S17} B. Simon, {\sc Unitaries permuting two orthogonal projections}, Linear Algebra
Appl. 528 (2017), 436-441.

\bibitem{WDD09} Y. Wang, H. Du, Y. Dou, {\sc On the index of Fredholm pairs of idempotents}, Acta Math. Sin. (Engl. Ser.) 25 (2009), 679-686.

\bibitem{W86} P.Y. Wu, {\sc Approximation by partial isometries}, Proc. Edinb. Math. Soc. 29 (1986),
 255-261.
\end{thebibliography}
\end{document}